\documentclass{article}

\usepackage{amssymb, amsmath}
\usepackage{color}
\usepackage{amsmath,stackengine}

\definecolor{red}{rgb}{0.9,0,0}

\definecolor{azul}{rgb}{0,0,1}

\usepackage{indentfirst}
\usepackage{graphicx}
\newtheorem{theorem}{Theorem}[section]
\newtheorem{proposition}{Proposition}[section]

\newtheorem{lemma}{Lemma}[section]

\usepackage[left=2.9cm,right=2.9cm,top=2.8cm,bottom=2.8cm]{geometry}

\newcommand{\fim}{\hfill $\Box$}

\def\dem{\noindent{\it Proof. }}

\def\fd{\hfill{\vbox to 7pt{\hbox to 7pt{\vrule height 7pt width 7pt}}}}

\def\:{\colon}

\title{On a quasilinear elliptic problem involving the 1-laplacian operator and a discontinuous nonlinearity}

\author{Marcos T. O. Pimenta
\\{\footnotesize E-mail address: marcos.pimenta@unesp.br}
\\{\footnotesize Departamento de Matem\'atica e Computa\c{c}\~ao, Universidade Estadual Paulista - Unesp}
\\ {\footnotesize CEP: 19060-900, Presidente Prudente - SP, Brazil},
\\Gelson Concei\c{c}\~ao G. dos Santos
\\{\footnotesize E-mail address: gelsonsantos@ufpa.br}
\\{\footnotesize Faculdade de Matem\'atica, Universidade Federal do Par\'a,
}
\\ {\footnotesize CEP:  66075-110, Bel\'em - PA, Brazil.}
\\Jo\~ao R. dos Santos J\'unior
\\{\footnotesize E-mail address: joaojunior@ufpa.br}
\\{\footnotesize Faculdade de Matem\'atica, Universidade Federal do Par\'a,
}
\\ {\footnotesize CEP:  66075-110, Bel\'em - PA, Brazil.}
}

\date{}

\begin{document}

\pretolerance10000

\maketitle

\begin{abstract}
In this work, we study a quasilinear elliptic problem involving the 1-laplacian operator, with a discontinuous, superlinear and subcritical nonlinearity involving the Heaviside function $H(\cdot - \beta)$. Our approach is based on an analysis of the associated p-laplacian problem, followed by a thorough analysis of the asymptotic behaviour or such solutions as $p \to 1^+$. We study also the asymptotic behaviour of the solutions, as $\beta \to 0^+$ and we prove that it converges to a solution of the original problem, without the discontinuity in the nonlinearity.
\end{abstract}



\numberwithin{equation}{section}
\bibliographystyle{plain}
\maketitle

\section{Introduction}
\label{Introduction}

In this work we study the following quasilinear elliptic equation
\begin{equation}
\label{P1intro}
\begin{cases}
\displaystyle -\Delta_1 u = H(u - \beta)|u|^{q-2}u &\mbox{in} \quad \mbox{$\Omega$} \\
\qquad u =0 &\mbox{on} \quad \mbox{${\partial}\Omega$,}
\end{cases}
\end{equation}
where the 1-laplacian operator is formally defined as $\Delta_1 u = \mbox{div}\left(Du/|Du|\right)$, $\Omega \subset \mathbb{R}^N$ is a bounded domain with a Lipschitz boundary, $N \geq 2$, $1 < q < N/(N-1),$
$\beta>0$ is a real parameter
and $H: \mathbb{R}\rightarrow\mathbb{R}$ is the Heaviside function $H(t) = 1$ if $t \geq 0$ and $H(t) = 0$ otherwise.

%The presence of the function $H$ in \eqref{P1intro} characterizes it as a problem with discontinuous nonlinearity

In recent decades, the study of nonlinear partial differential equations with discontinuous nonlinearities has attracted the attention of several researchers. One of the reasons to study such equations is due to many free boundary problems arising in mathematical physics which can be stated in this form. Among these problems, we have the obstacle problem, the seepage surface problem, and the Elenbaas equation, see \cite{chang0,chang,Chang1}. For more applications see \cite{Ambrosetti-Dobarro}.
Several techniques have been developed or applied to study this kind of problem, such as variational methods for nondifferentiable functionals, lower and upper solutions, dual variational principle, global branching, Palais principle of symmetric criticality for locally Lipschitz functional and the theory of multivalued mappings. See for instance,
Alves, Yuan and Huang \cite{Alves-Huang}, Alves, Santos and Nemer \cite{Alves-Nemer}, Ambrosetti and Badiale \cite{Ambrosetti-Badiale}, Ambrosetti, Calahorrano and Dobarro \cite{Ambrosetti-Dobarro}, Ambrosetti and Turner \cite{ambrosettiturner}, Anmin and Chang \cite{Anmin}, Arcoya and Calahorrano \cite{Arcoya-Calahorrano},  Cerami \cite{cerami}, Chang \cite{chang0,chang,Chang1}, Clarke \cite{Clarke,Clarke1}, Gazzola and R\v{a}dulescu \cite{Radulescu}, Krawcewicz and Marzantowicz \cite{Krawcewicz}, Molica Bisci and Repov\v{s}
\cite{molica-repov}, R\v{a}dulescu \cite{Radulescu2}, dos Santos and Figueiredo \cite{gelsongiovany} and their references.

As far as problems involving the 1-laplacian operator are concerned, there are at least two approaches one can follow. The first one is based on the study of the energy functional associated to the problem, which is defined in $BV(\Omega)$, whenever one can write it as the difference of a convex and locally Lipschitz functional and a $C^1$ one. Then, one can use the tools of nonsmooth nonlinear analysis (see \cite{chang,Clarke1,Szulkin}) to find critical points of such energy functional. Note that, in studying \eqref{P1intro}, this is far from being an option for us, since the energy functional associated to \eqref{P1intro} would be defined in $BV(\Omega)$, and given by
$$
I_H(u) = \|u\|_{BV(\Omega)} - \mathcal{F}_\beta(u),
$$
where $\mathcal{F}_\beta(u)=\int_{\Omega} F_H(u)dx$, with $f_H(s)=H(s-\beta)|s|^{q-2}s$ and $F_H(t) = \int_0^t f_H(s)dx$.
 Hence, since $\mathcal{F}_\beta$ is not a $C^1$ functional defined on $BV(\Omega)$, it would be tricky to show that a critical point of $I_H$ satisfies \eqref{P1intro} in some sense, since in this case, we could not use variational inequalities to follow the standard approach, which is based on that one proposed by \cite{Szulkin}.

 Fortunately, there is another approach which is based on the study of \eqref{P1intro}, with the 1-laplacian substituted by the p-laplacian operator, for $p > 1$. Then, one can use standard arguments to solve the associated problem and then studying the family of such solutions as $p \to 1^+$. To the best of our knowledge, the pioneering works involving this operator were written by F. Andreu, C. Ballesteler, V. Caselles and J.M. Maz\'on in a series of papers (among them \cite{AndreuBallesterCasellesMazonJFA,AndreuBallesterCasellesMazon,AndreuBallesterCasellesMazonDIE}), which gave rise to the monograph \cite{AndreuCasellesMazonMonograph}. Among the very first works on this issue we should also cite the works of Kawohl \cite{Kawohl} and Demengel \cite{Demengel}.

Before to state our main result, let us define what we mean by a solution of the
problem \eqref{P1intro}.
Inspired by locally Lipschitz continuous functionals \cite{chang,Clarke,Clarke1,grossinho}
and
Anzellotti-Frid-Chen's Pairing Theory \cite{Anzellotti1983,FridChen} (see Subsections \ref{Liploc} and \ref{BV} for more details),
we say that $u \in BV(\Omega)$ is a bounded variation solution of \eqref{P1intro}, if there exist
$ \rho \in L^{\frac{q}{q-1}}(\Omega) $
and ${\textbf z} \in X_N(\Omega)$ with $\|{\textbf z}\|_\infty \leq 1$,
such that
\begin{equation}
\left\{
\begin{array}{rcl}
- \mbox{div}\, {\textbf z} & = & \rho \quad  \mbox{ in $\mathcal{D}'(\Omega)$,}\\
({\textbf z},Du) & = & |Du| \quad \mbox{in $\mathcal{M}(\Omega)$},\\
{[{\textbf z},\nu]} & \in & \mbox{$sign(-u)$ \quad $\mathcal{H}^{N-1}$-a.e. on $\partial \Omega,$}
\end{array}\label{sol.u}
 \right.
\end{equation}
and it holds that, for almost every $x\in \Omega,$
\begin{equation}
\rho(x) \in
\begin{cases}
\{0\}, & \ \mbox{if} \ u(x)<\beta,\\
 [0, \beta^{q-1}], & \ \mbox{if} \ u(x)=\beta\\
 \{u(x)^{q-1}\}, & \ \mbox{if} \ u(x)>\beta.
 \end{cases}\label{sol.rho}
\end{equation}
It is important to point out here that if the set $\{x\in\Omega: u(x)=\beta\}$ has zero Lebesgue measure, then $\rho(x)=H(u(x)-\beta)|u(x)|^{q-2}u(x)$ for almost every $x\in \Omega$.

Motivated by the works previously mentioned, our first main result is the following.

%%%%%%%%%%%%%%%%%%%%%%%%%%%%%%%%%%%%%%%%%%%%%%%%%%%%%%%%%%%%%%%%%%%%%%%%%%%%%%%%%%%%%%%%%%%%%%%%%%%%%%%%%%%%%%%%%%%%%%%%%%%%%%%%%%%%%%%%%%%%%%%%%%%%%%%%%%%%%%%%%%%%%%%%%%%%%%%%%%%%%%%
%                                                                         MAIN RESULTS
%%%%%%%%%%%%%%%%%%%%%%%%%%%%%%%%%%%%%%%%%%%%%%%%%%%%%%%%%%%%%%%%%%%%%%%%%%%%%%%%%%%%%%%%%%%%%%%%%%%%%%%%%%%%%%%%%%%%%%%%%%%%%%%%%%%%%%%%%%%%%%%%%%%%%%%%%%%%%%%%%%%%%%%%%%%%%%%%%%%%%%%

\begin{theorem}
\label{existenceP1}
Suppose that $N \geq 2$ and $1 < q < N/(N-1)$. Then, for each $\beta > 0,$ \eqref{P1intro} admits at least one nonnegative and nontrivial solution
$u_\beta \in BV(\Omega)\cap L^{\infty}(\Omega)$,
in the sense of \eqref{sol.u}.
\end{theorem}

In the scope of the last theorem, a question which naturally arises is about the behavior of the solutions $u_\beta$, as $\beta \to 0^+$. In fact, one should expect that $u_\beta$ converges in some sense, as $\beta \to 0^+$, to a solution of the following problem
\begin{equation}
\label{P2intro}
\begin{cases}
\displaystyle -\Delta_1 u = |u|^{q-2}u &\mbox{in} \quad \mbox{$\Omega$} \\
\qquad u =0 &\mbox{on} \quad \mbox{${\partial}\Omega$.}
\end{cases}
\end{equation}
In the next theorem, we prove that this in fact occurs.

\begin{theorem}
\label{theorem2}
For each $\beta > 0$, let $u_\beta$ be the solution given in Theorem \ref{existenceP1}. Then there exists a nontrivial and nonnegative solution of \eqref{P2intro}, $u_0 \in BV(\Omega)$, such that, as $\beta \to 0^+$,
$$
u_\beta \to u_0 \quad \mbox{in $L^r(\Omega)$, for all $1 \leq r < N/(N-1)$ and also a.e. in $\Omega$.}
$$
Moreover, there exist positive constants $\mu$ and $\beta_{0},$ such that
\begin{equation}\label{meas}
|\{x\in \Omega: u_\beta(x)>\beta\}|\geq
\mu, \ \mbox{ for all } \ \beta\in(0, \beta_{0}),
\end{equation}
where $|A|$ denotes the measure of a measurable set $A\subset \mathbb{R}^N$.
\end{theorem}

Note that the last part of the theorem guarantees that the set $\{u_\beta > \beta\}$ does not shrink as $\beta \to 0$, that is, $\|u_\beta\|_{L^\infty(\Omega)}> \beta$, for $\beta$ small enough. Such an information is quite relevante because it ensures that, at least for $\beta$ small, $u_{\beta}$ is in fact a solution of a problem involving a discontinuous nonlinearity.

%%%%%%%%%%%%%%%%%%%%%%%%%%%%%%%%%%%%%%%%%%%%%%%%%%%%%%%%%%%%%%%%%%%%%%%%%%%%%%%%%%%%%%%%%%%%%%%%%%%%%%%%%%%%%%%%%%%%%%%%%%%%%%%%%%%%%%%%%%%%%%%%%%%%%%%%%%%%%%%%%%%%%%%%%%%%%%%%%%%%%%%
%                                                                            OUTLINE
%%%%%%%%%%%%%%%%%%%%%%%%%%%%%%%%%%%%%%%%%%%%%%%%%%%%%%%%%%%%%%%%%%%%%%%%%%%%%%%%%%%%%%%%%%%%%%%%%%%%%%%%%%%%%%%%%%%%%%%%%%%%%%%%%%%%%%%%%%%%%%%%%%%%%%%%%%%%%%%%%%%%%%%%%%%%%%%%%%%%%%%

The existence of positive solution for 
\eqref{P1intro}  with $\beta=0$ (i.e., \eqref{P2intro}) was recently studied by Molino-Segura in \cite{SeguraMolino}.
Due to the discontinuity in \eqref{P1intro}, caused by the Heaviside function (with $\beta>0$), we cannot use the classical critical point theory for $C^1$ functionals as in \cite{SeguraMolino}.
For this reason, motivated by
\cite{chang,ambrosettiturner,Arcoya-Calahorrano,Clarke,Clarke1,Chang1-Lapla}, we combine variational methods for nondifferentiable functionals
with the approximation argument of \cite{SeguraMolino}.

In Theorem \ref{existenceP1}, to prove the boundedness of the solutions, we use Moser's iteration method (see \cite{Moser}) and
a careful analysis of some constants
to obtain a uniform estimate in the $L^\infty(\Omega)-$norm of the solutions of the approximate problem.
These estimates were essential in our arguments to ensure that the solution of problem \eqref{P1intro} is nontrivial.

This paper is organized as follows. In Section \ref{Sectionpreliminaries} we present some definitions and basic results about functions of bounded variation and the nonlinear analysis involving
nonsmooth functionals. In Sections \ref{SectionExistenceP1} and \ref{Prooftheorem2}, we present the proofs of Theorem \ref{existenceP1} and \ref{theorem2}, respectively.

\section{Preliminaries}
\label{Sectionpreliminaries}

\subsection{Main properties of $BV(\Omega)$ space}\label{BV}

First of all let us introduce the space of functions of bounded variation, $BV(\Omega)$, where $\Omega \subset \mathbb{R}^N$ is a domain. We say that $u \in BV(\Omega)$, or is a function of bounded variation, if $u \in L^1(\Omega)$, and its distributional derivative $Du$ is a vectorial Radon measure, i.e.,
$$
BV(\Omega) = \left\{u \in L^1(\Omega); \, Du \in \mathcal{M}(\Omega,\mathbb{R}^N)\right\}.
$$
It can be proved that $u \in BV(\Omega)$ if and only if $u \in L^1(\Omega)$ and
$$
\int_{\Omega} |Du| := \sup\left\{\int_{\Omega} u \, \mbox{div}\, \phi dx; \, \, \phi \in C^1_c(\Omega,\mathbb{R}^N), \, \|\phi\|_\infty \leq 1\right\} < +\infty.
$$

The space $BV(\Omega)$ is a Banach space when endowed with the norm
$$\|u\|_{BV} := \int_{\Omega} |Du| + \int_{\Omega}|u|dx,$$
which is continuously embedded into $L^r(\Omega)$ for all $\displaystyle r \in \left[1,1^*\right]$, where $1^* = N/(N-1)$. Since the domain $\Omega$ is bounded, it holds also the compactness of the embeddings of $BV(\Omega)$ into $L^r(\Omega)$ for all $\displaystyle r \in \left[1,1^*\right)$.

The space $C^\infty(\overline{\Omega})$ is not dense in $BV(\Omega)$ with respect to the strong convergence. However, with respect to the {\it strict convergence}, it does. We say that $(u_n) \subset BV(\Omega)$ converges to $u \in BV(\Omega)$ in the sense of the strict convergence, if
$$
u_n \to u, \quad \mbox{in $L^1(\Omega)$}
$$
and
$$
\int_{\Omega}|Du_n| \to \int_{\Omega}|Du|,
$$
as $n \to \infty$. In \cite{AmbrosioFuscoPallara} one can see also that it is well defined a trace operator $BV(\Omega) \hookrightarrow L^1(\partial \Omega)$, in such a way that
$$
\|u\| := \int_{\Omega} |Du| + \int_{\partial\Omega}|u|d\mathcal{H}^{N-1},
$$
is a norm equivalent to $\|\cdot\|_{BV}$.

Given $u \in BV(\Omega)$, we can decompose its distributional derivative as
$$
Du = D^au + D^su,
$$
where $D^au$ is absolutely continuous with respect to the Lebesgue measure $\mathcal{L}^N$, while $D^su$ is singular with respect to the same measure. Moreover, we denote de total variation of $Du$, as $|Du|$.

In several arguments we use in this work, it is mandatory to have a sort of {\it Green's formula} to expressions like $w \, \mbox{div}({\textbf z})$, where ${\textbf z} \in L^\infty(\Omega,\mathbb{R}^N)$, $\mbox{div}({\textbf z}) \in L^N(\Omega)$ and $w \in BV(\Omega)$. For this we have to somehow deal with the {\it product} between ${\textbf z}$ and $Dw$, which we denote by $({\textbf z},Dw)$. This can be done through the {\it pairings theory}, developed by Anzellotti in \cite{Anzellotti1983} and independently by Frid and Chen in \cite{FridChen}. Below, we describe the main results of this theory.

Let us denote
$$
X_N(\Omega) = \left \{{\textbf z} \in L^\infty(\Omega,\mathbb{R}^N); \, \mbox{div}({\textbf z}) \in L^N(\Omega) \right \}.
$$
For ${\textbf z} \in X_N(\Omega)$ and $w \in BV(\Omega)$, we define the distribution $({\textbf z},Dw) \in \mathcal{D}'(\Omega)$ as
$$
\langle ({\textbf z},Dw), \varphi \rangle := - \int_\Omega w \varphi \, \mbox{div}({\textbf z})dx - \int_\Omega w {\textbf z}\cdot \nabla \varphi dx,
$$
for every $\varphi \in \mathcal{D}(\Omega)$. With this definition, it can be proved that $({\textbf z},Dw)$ is in fact a Radon measure such that
\begin{equation}
\left|\int_B ({\textbf z},Dw)\right| \leq \|{\textbf z}\|_\infty \int_B |Dw|,
\label{zDuinequality}
\end{equation}
for every Borel set $B \subset \Omega$.

In order to define an analogue of the Green's Formula, it is also necessary to describe a weak trace theory for ${\textbf z}$. In fact, there exists a trace operator $\left[ \cdot, \nu\right]: X_N(\Omega) \to L^\infty(\partial \Omega)$ such that
\begin{equation}
\|\left[{\textbf z}, \nu\right]\|_{L^\infty(\partial \Omega)} \leq \|{\textbf z}\|_\infty
\label{ineqznu}
\end{equation}
and, if ${\textbf z} \in C^1(\overline{\Omega}_\delta,\mathbb{R}^N)$,
$$
\left[{\textbf z}, \nu\right](x) = {\textbf z}(x)\cdot \nu(x) \quad \mbox{on $\Omega_\delta$,}
$$
where by $\Omega_\delta$ we denote a $\delta$-neighborhood of $\partial \Omega$. With these definitions, it can be proved that the following Green's formula holds for every ${\textbf z} \in X_N(\Omega)$ and $w \in BV(\Omega)$,
\begin{equation}
\int_\Omega w\, \mbox{div}({\textbf z})dx + \int_\Omega ({\textbf z},Dw) = \int_{\partial \Omega}[{\textbf z},\nu] w d\mathcal{H}^{N-1}.
\label{GreenFormula}
\end{equation}

\subsection{Nonlinear analysis on non-differentiable functionals}\label{Liploc}

In this subsection, for the reader's convenience,  we recall some
definitions and basic results on the critical point theory of
locally Lipschitz continuous functionals (that is based on the subdifferential theory
of Clarke \cite{Clarke,Clarke1})
as developed by Chang
\cite{chang}, Clarke \cite{Clarke,Clarke1} and Grossinho and Tersian
\cite{grossinho}.

Let $E$ be a real Banach space. A functional $I:E \rightarrow {\mathbb{R}}$ is locally Lipschitz continuous, $I \in
Lip_{loc}(E, {\mathbb{R}})$ for short, if given $u \in E$ there is
an open neighborhood $V := V_u \subset E$ and some constant  $M =
M_V > 0$ such that
$$
\mid I(v_2) - I(v_1) \mid \leq M \| v_2-v_1 \|, \ v_i
\in V, \ i = 1,2.
$$
The directional derivative of $I$ at $u$ in the direction of $v
\in E$ is defined by
$$
I^0(u;v)=\displaystyle \displaystyle \limsup_{h \to 0,~\sigma
\downarrow 0} \frac{I(u+h+ \sigma v)-I(u+h)}{\sigma}.
$$
Hence  $I^0(u;.)$ is continuous, convex and its subdifferential
at $z\in E$ is given by
$$
\partial I^0(u;z):=\big\{\mu\in E^*;I^0(u;v)\geq I^0(u;z)+
\langle\mu,v-z\rangle, \ v\in E\big\},
$$
where $\langle.,.\rangle$ is the duality pairing between $E^*$
and $E$. The generalized gradient of $I$ at $u$ is the set
$$
\partial I(u):=\big\{\mu\in E^*; \langle \mu,v\rangle\leq I^0(u;v), \ v\in
E\big \}.
$$
Since $I^0(u;0) = 0$, $\partial I(u)$ is the subdifferential of
$I^0(u;.)$ in $0$.

It is also known that $\partial I(u)\subset E^*$ is convex, non-empty and weak*-compact
and it is well defined
\begin{eqnarray}\label{mI}
\Lambda_I (u): = \min \big\{\parallel\mu\parallel_{E^*};\mu \in
\partial I(u) \big \}.
\end{eqnarray}

A critical point of $I$ is an element $u_\beta \in E$ such that $0\in
\partial I(u_\beta)$ and a critical value of $I$ is a real number $c$
such that $I(u_\beta)=c$ for some critical point $u_\beta \in E$.

We say that $I\in Lip_{loc}(E, \mathbb{R})$ satisfies the nonsmooth Palais-Smale condition at level $c\in\mathbb{R}$ (nonsmooth $(PS)_c$-condition for short), if the following holds: every sequence $(u_n)\subset E$, such that $I(u_n)\rightarrow c$ and $\Lambda_I(u_n)\rightarrow0$ has a strongly convergent subsequence.

\begin{proposition}\label{propri} (See \cite{Clarke,Clarke1,grossinho})
Let $I_1, I_2: E\rightarrow \mathbb{R}$ be locally Lipschitz functions, then:

\noindent{\bf $(i)$} $I_1+ I_2\in Lip_{loc}(E,\mathbb{R})$ and $ \partial(I_1+ I_2)(u) \subseteq \partial I_1(u) + \partial I_2(u),$ for all $u\in E.$

\noindent{\bf $(ii)$} $\partial (\lambda I_1)(u) = \lambda \partial I_1(u)$ for each $\lambda\in \mathbb{R}, u\in E.$

\noindent{\bf $(iii)$} Suppose that for each point $v$ in a neighborhood of $u$, $I_1$ admits a Gateaux
derivative $I_1'(v)$ and that $I_1' :E\rightarrow E^*$ is continuous, then $\partial I_1(u) = \{I_1'(u)\}.$

\end{proposition}

\begin{theorem}\label{TPM} (See \cite{Clarke,Clarke1,grossinho}) Let $E$ be a Banach space and let $I\in Lip_{loc}(E,\mathbb{R})$ with $I(0) = 0.$ Suppose there are
numbers $\alpha,r>0$ and $e\in E,$ such that

\noindent{\bf $(i)$} $I(u)\geq \alpha,$  for all $u\in E; \|u\|=r,$

\noindent{\bf $(ii)$} $I(e)<0$ and $\|e\|>r.$

\noindent Let
\begin{eqnarray}\label{level-minimax}
c=\inf_{\gamma\in\Gamma}\max_{t\in[0,1]} I(\gamma(t)) \mbox{ and } \Gamma=\{\gamma\in C([0,1],E): \gamma(0)=0\ \mbox{ and } \gamma(1)=e\}.
\end{eqnarray}
Then $c\geq\alpha$ and there is a sequence $(u_{n})\subset E$ satisfying
$$I(u_{n})\rightarrow c\ \mbox{ and } \ \Lambda_{I}(u_{n})\rightarrow0.$$
If, in addition, $I$ satisfies the
nonsmooth $(PS)_c$-condition, then $c$ is a critical value of $I.$
\end{theorem}

\medskip

%%%%%%%%%%%%%%%%%%%%%%%%%%%%%%%%%%%%%%%%%%%%%%%%%%%%%%%%%%%%%%%%%%%%%%%%%%%%%%%%%%%%%%%%%%%%%%%%%%%%%%%%%%%%%%%%%%%%%%%%%%%%%%%%%%%%%%%%%%%%%%%%%%%%%%%%%%%%%%%%%%%%%%%%%%%%%%%%%%%%%%%
%                                                                     AN AUXILIARY PROBLEM
%%%%%%%%%%%%%%%%%%%%%%%%%%%%%%%%%%%%%%%%%%%%%%%%%%%%%%%%%%%%%%%%%%%%%%%%%%%%%%%%%%%%%%%%%%%%%%%%%%%%%%%%%%%%%%%%%%%%%%%%%%%%%%%%%%%%%%%%%%%%%%%%%%%%%%%%%%%%%%%%%%%%%%%%%%%%%%%%%%%%%%%

\section{Proof of Theorem \ref{existenceP1}}
\label{SectionExistenceP1}

In this section, to prove our main result, we will consider a family of auxiliary problems involving the p-Laplacian operator and discontinuous nonlinearity.
We will use an approximation technique and variational methods for nondifferentiable functionals
inspired by Molino-Segura de Le\'on \cite{SeguraMolino}, Anzellotti-Frid-Chen \cite{Anzellotti1983,FridChen},
Arcoya-Calahorrano \cite{Arcoya-Calahorrano}, Ambrosetti-Turner \cite{ambrosettiturner}, Clarke \cite{Clarke} and
Chang \cite{chang}.

In order to get such solutions of \eqref{P1intro}, the first step is to consider the problem
\begin{equation}
\label{Pp}
\begin{cases}
\displaystyle -\mbox{div}\left(|\nabla u|^{p-2}\nabla u\right) = H(u - \beta)|u|^{q-2}{u} &\mbox{in} \quad \mbox{$\Omega$,} \\
\qquad u =0 &\mbox{on} \quad \mbox{${\partial}\Omega$.}
\end{cases}
\end{equation}

We say that $u_{p,\beta} \in W^{1,p}_0(\Omega)$ is a weak solution of \eqref{Pp}, if there exists $\rho_{p,\beta} \in L^{\frac{q}{q-1}}(\Omega)$, such that
\begin{equation}
\int_\Omega |\nabla u_{p,\beta}|^{p-2}\nabla u_{p,\beta} \nabla \varphi dx = \int_\Omega \rho_{p,\beta} \varphi dx, \ \mbox{ for all } \ \varphi \in W^{1,p}_0(\Omega),
\label{Ppweak}
\end{equation}
and it holds that, for almost every $x\in \Omega,$
\begin{equation}
\rho_{p,\beta}(x) \in
\begin{cases}
\{0\}, & \ \mbox{if} \ u_{p,\beta}(x)<\beta,\\
 [0, \beta^{q-1}], & \ \mbox{if} \ u_{p,\beta}(x)=\beta, \\
 \{u_{p,\beta}(x)^{q-1}\}, & \ \mbox{if} \ u_{p,\beta}(x)>\beta.
 \end{cases}\label{rhop}
\end{equation}

%Using the nonsmooth critical point theory, it is possible to prove that the problem \eqref{Pp} has at least one solution $u_{p,\beta}\in W_0^{1,p}(\Omega),$ which can be obtained by the nonsmooth version of the Mountain Pass Theorem, see Theorem \ref{TPM}. Since we are interested in some properties of these solutions, added to the for the reader's convenience, we will give the details of the proof.

Inspired by Arcoya and Calahorrano \cite{Arcoya-Calahorrano}, which proved the existence of solution for a sublinear version of \eqref{Pp} (see also Ambrosetti and Turner \cite{ambrosettiturner}), we will use the nonsmooth critical point theory to prove that problem \eqref{Pp} has at least one solution $u_{p,\beta}\in W_0^{1,p}(\Omega)$, which will be obtained by the nonsmooth version of the Mountain pass theorem (see Theorem \ref{TPM}).
Furthermore, we will prove some properties of this solution that will be useful to prove the existence of a solution to problem \eqref{P1intro}.
To achieve this goal, first note that by Chang's results \cite[Theorem 2.1 and Theorem 2.3]{chang}, %and \cite[Example 1, pg 107]{chang}, it follows that
the functional $\mathcal{F}_\beta: L^q(\Omega)\rightarrow \mathbb{R}$
given by
$$\mathcal{F}_\beta(u)=\int_{\Omega} F_\beta(u)dx, \mbox{ with $f_\beta(s) = H(s-\beta)|s|^{q-2}s$ and $F_\beta(t) = \int_0^t f_\beta(s)dx$,}
$$
is locally Lipschitz and
\begin{equation}\label{gradF}
\partial \mathcal{F}_\beta(u)=[\underline{f}_\beta(u),\overline{f}_\beta(u)] \ \mbox{ a.e. in } \Omega,
\end{equation}
where
$$\underline{f}_\beta(t)=
\lim_{r\rightarrow0^+} \mbox{ess inf}\{f_\beta(s): |t-s| < r\} \ \mbox{ and } \overline{f}_\beta(t) = \lim_{r\rightarrow0^+} \mbox{ess sup}\{f_\beta(s): |t-s| < r\}.$$
It is clear that
\begin{equation}
[\underline{f}_\beta(t),\overline{f}_\beta(t)]=\begin{cases}
\{0\}, & \ \mbox{if} \ t<\beta,\\
 [0, \beta^{q-1}], & \ \mbox{if} \ t=\beta,\\
 \{t^{q-1}\}, & \ \mbox{if} \ t>\beta.
 \end{cases}\label{inclu}
\end{equation}

The associated functional for \eqref{Pp} is $J_{p,\beta}: W_{0}^{1,p}(\Omega)\rightarrow \mathbb{R}$, given by
\begin{equation}\label{JpH}
J_{p,\beta}(u)=Q_p(u)- \mathcal{F}_\beta\bigl|_{W_0^{1,p}}(u), \ \mbox{where} \ Q_p(u)=\frac{1}{p}\int_{\Omega}|\nabla u|^pdx.
\end{equation}

Due to the presence of the Heaviside function $H,$ the functional $J_{p,\beta}$ is not Fr\'echet differentiable, but is locally Lipschitz on $W_0^{1,p}(\Omega)$. Moreover, by \cite[Theorem 2.2]{chang} we have $\partial \big(\mathcal{F}_{H}\bigl|_{W_0^{1,p}}\big)(u)=\partial\mathcal{F}_\beta(u),$ for all $u\in W_0^{1,p}(\Omega).$ Hence,
by
Proposition \ref{propri},
\begin{equation}\label{grad}
\partial J_{p,\beta}(u)=\{Q_p'(u)\}-\partial \mathcal{F}_\beta(u) \ \mbox{ for all } u\in W_0^{1,p}(\Omega),
\end{equation}
and therefore, by \eqref{gradF}, \eqref{inclu} and \eqref{grad},  critical points of $J_{p,\beta},$ in the sense of the nonsmooth critical point theory, will give rise to
solutions of \eqref{Pp}.

\medskip

Since we want to find a nontrivial solution of \eqref{P1intro} by using the solutions $u_{p,\beta}$ of \eqref{Pp}
by passing to the limit as $p\rightarrow 1^{+},$
in what follows, we will consider $p\in(1,\overline{p})$ for some $\overline{p}\in(1,q)$ fixed.

%We remark that it would be natural to search for solutions of \eqref{P1intro} by studying the existence of critical points of J that are in fact solutions of (P).  But there are some technical difficulties in applying the critical point theory

\begin{lemma}\label{geometria} For each $p\in(1,\overline{p})$ and $\beta>0,$ the functional $J_{p,\beta}$
satisfies the geometric conditions of the Mountain pass theorem. More precisely,

\noindent{\bf $(i)$} There exist $r,\alpha>0$, which are independent of $\beta$, such that $J_{p,\beta}(u)\geq \alpha$ for all $u\in W_0^{1,p}(\Omega)$ with $\|u\|_{W_0^{1,p}(\Omega)}=r.$ Moreover, $\alpha$ can be chosen also independent of $p$.

\noindent{\bf $(ii)$} There exists $e = e(\beta)\in C_0^{\infty}(\Omega) $ such that $J_{p,\beta}(e) < 0$ and $\|e\|_{W_0^{1,p}(\Omega)}> r.$
\end{lemma}
\dem

By  H{\"o}lder's inequality,
$$
J_{p,\beta}(u) \geq \frac{1}{p}\|\nabla u\|_{L^p(\Omega)}^p - \biggl(\int_{\Omega}|u|^{p^*}\biggl)^{\frac{q}{p^*}}|\Omega|^{\frac{p^*-q}{p^*}},  \mbox{ for all } u\in W_0^{1,p}(\Omega).
$$
Since, by \cite[Proof of Theorem 7.10]{Gilbarg}, for each $u\in W_0^{1,p}(\Omega),$
\begin{eqnarray}\label{p*}
\|u\|_{L^{p^*}(\Omega)}\leq \frac{\theta}{\sqrt{N}}\|\nabla u\|_{L^{p}(\Omega)}, \mbox{ where } \theta=\frac{p(N-1)}{N-p},
\end{eqnarray}
we have,
\begin{eqnarray*}
J_{p,\beta}(u) &\geq &\frac{1}{p}\|\nabla u\|_{L^{p}(\Omega)}^p -
|\Omega|^{\frac{p^*-q}{p^*}}\biggl(\frac{p(N-1)}{\sqrt{N}(N-p)}\biggl)^{q}
\|\nabla u\|_{L^{p}(\Omega)}^q\\
&\geq & \frac{1}{\overline{p}}\|\nabla u\|_{L^{p}(\Omega)}^p -
C\|\nabla u\|_{L^{p}(\Omega)}^q,
\mbox{ for all } u\in W_0^{1,p}(\Omega),\\
\end{eqnarray*}
where $C=\bigl(\overline{p}(N-1)/\sqrt{N}(N-\overline{p})\bigl)^{q}\max\{1,|\Omega|\}.$

Note that
$$
\frac{r^p}{\overline{p}} - Cr^{q}\geq \frac{r^q}{\overline{p}} \ \mbox{ if and only if  }
0<r\leq\biggl(\frac{1}{\overline{p}C+1}\biggl)^{\frac{1}{q-p}}.
$$
Then, by choosing $r=\biggl(\frac{1}{\overline{p}C+1}\biggl)^{\frac{1}{q-\overline{p}}}$ and $\alpha=r^{q}/\overline{p},$
we conclude that $(i)$ holds.

Now, let $\varphi\in C^{\infty}_0(\Omega)$ be such that $|\{\varphi>\beta\}|>0,$
where $\{\varphi>\beta\}$ denotes the
set  $\{x\in \Omega : \varphi(x) > \beta\}.$
For each $t\geq1,$ we get
\begin{eqnarray*}
\begin{aligned}
J_{p,\beta}(t\varphi) & = \frac{t^p}{p}\| \varphi \|_{W_0^{1,p}(\Omega)}^p- \int_{\Omega}F_\beta(t\varphi)dx \\
& \leq  \frac{t^p}{p}\| \varphi \|_{W_0^{1,p}(\Omega)}^p-
\frac{t^q}{q}\int_{\{\varphi>\beta\}}\varphi^q  dx
+\frac{\beta^{q}}{q}|\Omega|,
\end{aligned}
\end{eqnarray*}
which implies in the existence of $e$ satisfying $(ii).$

\fim

\begin{lemma}\label{lema-PS} For each $p\in(1,\overline{p})$ and $\beta>0,$ $J_{p,\beta}$ satisfies the nonsmooth Palais-Smale condition.
\end{lemma}
\dem
Let $(u_n)\subset W_0^{1,p}(\Omega)$ be a $(PS)_c$ sequence for $J_{p,\beta},$ that is,
$J_{p,\beta}(u_n) \to c$  and $\Lambda_{J_{p,\beta}}(u_{n})\rightarrow 0,$ where $\Lambda_{J_{p,\beta}}$ is defined in \eqref{mI}.
Hence, it follows from \eqref{mI} and \eqref{grad} that there exists $(\mu_{n})\subset\partial J_{p,\beta}(u_{n})$ such that
$$
 \|\mu_{n}\|_{*}=\Lambda_{J_{p,\beta}}(u_{n})=o_{n}(1) \mbox{ and } \mu_{n}=Q'_p(u_{n})-\rho_{n},
$$
where  $\rho_{n}\in\partial\mathcal{F}_\beta(u_{n}).$ Then,
\begin{eqnarray}\label{bound}
\begin{aligned}
c + 1 + \| u_n\|_{W_0^{1,p}(\Omega)} & \geq  J_{p,\beta}(u_n) - \frac{1}{q}\left\langle \mu_n, u_n \right\rangle + o_n(1)\\
 & =  J_{p,\beta}(u_n) - \frac{1}{q}\left\langle Q_p'(u_n)-\rho_n, u_n \right\rangle + o_n(1)\\
 & =  \left(\frac{1}{p}-\frac{1}{q}\right)\| u_n\|^p_{W_0^{1,p}(\Omega)} +\int_{\Omega}\left(\frac{1}{q}\rho_n u_n - F_\beta(u_n)\right)  dx  + o_n(1).
 \end{aligned}
\end{eqnarray}

Moreover, note that by \eqref{rhop} and \eqref{gradF}, we have
\begin{equation}\label{>0}
 \int_{\Omega}\left(\frac{1}{q}\rho_n u_n - F_{\beta}(u_n)\right) dx = \frac{\beta}{q}\int_{\{u_n=\beta\}}\rho_n  dx + \frac{\beta^q}{q}|\{u_n > \beta\}| \geq 0.
\end{equation}
Hence,
\begin{eqnarray}
c + 1 + \| u_n\|_{W_0^{1,p}(\Omega)} \geq \left(\frac{1}{p}-\frac{1}{q}\right)\| u_n\|^p_{W_0^{1,p}(\Omega)} + o_n(1),
\end{eqnarray}
which implies that the sequence $(u_n)$ is bounded in $W_0^{1,p}(\Omega)$.
Thus, by Sobolev embedding theorems, passing to a subsequence if necessary, we obtain
\begin{equation}\label{ineq8}
\begin{cases}
u_n \rightharpoonup u \ \ \mbox{in} \ W_0^{1,p}(\Omega), \  u_n \rightarrow u \ \mbox{ in } \ L^{s}(\Omega), \\
u_n(x) \rightarrow u(x) \ \ \mbox{a.e in} \ \Omega, \\
|u_n(x)| \leq h(x) \ \mbox{ for some } \ h \in L^{s}(\Omega), s\in [1,p^*:=\frac{Np}{N-p}).
\end{cases}
\end{equation}

Using a similar argument than \cite[pg. 1071]{Arcoya-Calahorrano}, we conclude that $J_{p,\beta}$ satisfies the nonsmooth Palais-Smale condition.

\fim

Let us define the functional $I_{p,\beta}: W^{1,p}_0(\Omega) \to \mathbb{R}$, given by
$$
I_{p,\beta}(u) := J_{p,\beta}(u) + \frac{(p-1)}{p}|\Omega|.
$$

Note that, by Lemma \ref{geometria}, Lemma \ref{lema-PS} and Theorem \ref{TPM}, $I_{p,\beta}$ has a critical point  $u_{p,\beta} \in W_{0}^{1,p}(\Omega)$ at the level
$$
c_{p,\beta}=\inf_{\gamma\in\Gamma}\max_{t\in[0,1]} I_{p,\beta}(\gamma(t))\;\mbox{ with }\; \Gamma=\{\gamma\in C([0,1],W_0^{1,p}(\Omega)): \gamma(0)=0\;\mbox{and}\;\gamma(1)=e\},
$$
that is,
\begin{equation}\label{cp}
0\in \partial I_{p,\beta}(u_{p,\beta}) \mbox{  and } I_{p,\beta}(u_{p,\beta}) = c_{p,\beta}.
\end{equation}
 Hence, there exists $\rho_{p,\beta} \in L^{\frac{q}{q-1}}(\Omega) $ such that
$u_{p,\beta}$ and $\rho_{p,\beta}$ satisfy \eqref{Ppweak} and \eqref{rhop}.
Moreover,
testing \eqref{Ppweak} with $\varphi = u_{p,\beta}^-:=\min\{u_{p,\beta},0\}$
and using \eqref{rhop} we have $\|u^{-}_{p,\beta}\|^p_{W^{1,p}_0(\Omega)}=0,$ which implies that
$u_{p,\beta}(x)\geq0$ a.e. in $\Omega.$

\begin{lemma}\label{bounded} Let $u_{p,\beta}$ be given in \eqref{cp}. Then the family $(u_{p,\beta})_{1 < p < \overline{p}}$ is bounded in $BV(\Omega).$
\end{lemma}
\dem

By Young's inequality, we have
$$
\int_{\Omega}|\nabla u|^{p_1}dx\leq\frac{p_1}{p_2}\int_{\Omega}|\nabla u|^{p_2}dx+\frac{p_2-p_1}{p_2}|\Omega|,  \mbox{ for all } 1< p_1\leq p_2, u\in W_0^{1,p}(\Omega).
$$
Hence, $I_{p,\beta}$ is nondecreasing with respect to $p$ and
arguing as in \cite{SeguraMolino}, we conclude that
 $(I_{p,\beta}(u_{p,\beta}))_{1<p<\overline{p}}$ is increasing. Hence,
 \begin{equation}
 c_{p_1}\leq c_{p_2}
\label{ineqcp}
\end{equation}
 for all $1 < p_1\leq p_2.$
Note also that, by \eqref{cp},
\begin{eqnarray*}
c_{p,\beta}&=&I_{p,\beta}(u_{p,\beta})-\frac{1}{q}
\left\langle Q_p'(u_{p,\beta})-\rho_{p,\beta}, u_{p,\beta} \right\rangle \\
&=&\left(\frac{1}{p}-\frac{1}{q}\right)\int_{\Omega}| \nabla u_{p,\beta}|^pdx +\int_{\Omega}\left(\frac{1}{q}\rho_{p,\beta} u_{p,\beta} - F_\beta(u_{p,\beta})\right)dx
.
\end{eqnarray*}

From \eqref{>0} and \eqref{ineqcp}, it follows that
\begin{equation}
\int_{\Omega}|\nabla u_{p,\beta}|^p dx\leq C,  \mbox{ for all }  1<p<\overline{p},
\label{pnormestimate}
\end{equation}
where $C:=\frac{\overline{p}q}{q-\overline{p}}c_{\overline{p},\beta}>0$
is a constant independent of $p\in(1,\overline{p}).$

Applying once more Young's inequality, we obtain
\begin{eqnarray*}
\|u_{p,\beta}\| &\leq& \frac{1}{p}\int_\Omega|\nabla u_{p,\beta}|^{p}dx + \frac{p-1}{p}|\Omega|\\
&  \leq & \overline{C}+ |\Omega|,
\end{eqnarray*}
for some constant $\overline{C}>0,$ independent of $p.$

\fim

\begin{lemma}\label{regularity1}
For each $\beta>0,$ the function $u_{p,\beta}$ given in \eqref{cp} satisfies
\begin{eqnarray}
\|u_{p,\beta}\|_{L^\infty(\Omega)}\leq \overline{C},
\end{eqnarray}
for some constant $\overline{C}>0$ independent of $p\in(1,\overline{p}).$
\end{lemma}
\dem
Here to simplify the notation we put $u=u_{p,\beta}$ and $\rho_{p,\beta}=\rho.$
For obtain the $L^\infty$-estimate we will use the Moser's iteration \cite{Moser}
and a careful analysis of some constants.
For each $L>0$, we define
$$
u_L(x):= \left\{
\begin{array}{rcl}
u(x),& \mbox{if} & u(x)\leq L\\
L, & \mbox{if} & u(x)> L,
\end{array}
\right.
$$
$$z_{L,n}(x):=(u_L^{p(\gamma-1)}u)(x)\ \mbox{ and }\ w_L(x):=(uu_L^{\gamma-1})(x),$$
with $\gamma>1$ to be determined later. Choosing $\varphi=z_{L,n}$ in \eqref{Ppweak}, we get
\begin{eqnarray*}
\int_{\Omega}u_L^{p(\gamma-1)}|\nabla u|^{p}dx\displaystyle=-p(\gamma-1)\int_{\Omega}u_L^{p\gamma-p-1}u|\nabla u|^{p-2}\nabla u\nabla u_L dx+\int_{\Omega}\rho u u_L^{p(\gamma-1)}dx.
\end{eqnarray*}
Since
$$p(\gamma-1)\int_{\Omega}u_L^{p\gamma-p-1}u|\nabla u|^{p-2}\nabla u\nabla u_L dx=p(\gamma-1)\int_{\{u\leq L\}}u_L^{p(\gamma-1)}|\nabla u|^pdx\geq0$$
and $0\leq \rho(x)\leq  |u(x)|^{q-1}$ for almost every $x\in\Omega,$ see \eqref{rhop},
we obtain
\begin{eqnarray}\label{nLinea23}
\int_{\Omega}u_L^{p(\gamma-1)}|\nabla u|^{p}dx\leq \int_{\Omega}u^qu_L^{p(\gamma-1)}dx.
\end{eqnarray}
On the other hand, by \eqref{p*} it follows that
$$|w_L|_{L^{p^*}(\Omega)}^{p}\leq c_{p,\beta}\int_{\Omega}|\nabla w_L|^{p}dx
=c_{p,\beta}\int_{\Omega}
|\nabla(uu_L^{\gamma-1})|^{p}dx,
$$
where $c_{p,\beta}=\bigl(p(N-1)/(N-p)\bigl)^{p}.$
Thus,
$$
|w_L|_{L^{p^*}(\Omega)}^{p}\leq 2^pc_{p,\beta}\int_{\Omega}u_L^{p(\gamma-1)}|\nabla u|^pdx+2^pc_{p,\beta}(\gamma-1)^{p}\int_{\Omega}u_L^{p(\gamma-2)}u^p|\nabla u_L|^p dx,$$
hence, we get
\begin{eqnarray}\label{nLinea24}
|w_L|_{L^{p^*}(\Omega)}^p\leq 2^pc_{p,\beta}\gamma^p\int_{\Omega}u_L^{p(\gamma-1)}|\nabla u|^pdx.
\end{eqnarray}
Combining $(\ref{nLinea23})$ and $(\ref{nLinea24})$, we obtain
\begin{eqnarray*}
|w_L|_{L^{p^*}(\Omega)}^p\leq  2^pc_{p,\beta}\gamma^p\int_{\Omega}u^{q-p}(uu_L^{\gamma-1})^pdx,
\end{eqnarray*}
and so,
\begin{eqnarray*}
|w_L|_{L^{p^*}(\Omega)}^{p}\leq 2^pc_{p,\beta} \gamma^p\int_{\Omega}u^{q-p}w_L^pdx.
\end{eqnarray*}
Now we use the H{\"o}lder's inequality (with exponents $p^*/(q-p)$ and $p^*/(p^*-(q-p))$ to get that
\begin{eqnarray*}
|w_L|_{L^{p^*}(\Omega)}^p\leq 2^pc_{p,\beta}\gamma^p\Biggl(\int_{\Omega}u^{p^*}dx\Biggl)^{\frac{q-p}{p^*}}\Biggl(\int_{\Omega}w_L^{\frac{pp^*}{p^*-(q-p)}}dx\Biggl)^{\frac{p^*-(q-p)}{p^*}},
\end{eqnarray*}
where $p<\frac{pp^*}{p^*-(q-p)}<p^*.$
\\
The previous inequality, \eqref{p*} and \eqref{pnormestimate} imply that
\begin{eqnarray}\label{wL}
|w_L|_{L^{p^*}(\Omega)}^{p}\leq (2\theta_p)^p\gamma^p\Biggl(\int_{\Omega}w_L^{\alpha^*}dx\Biggl)^{\frac{p}{\alpha^*}},
\end{eqnarray}
where
\begin{eqnarray}\label{alphatheta}
\alpha^*:=\frac{pp^*}{p^*-(q-p)}\ \mbox{ and } \ \theta_p:=\biggl(p\frac{(N-1)}{(N-p)}\biggl)^{\frac{q}{p}}C^{\frac{q-p}{p}},
\end{eqnarray}
with the constant $C$ given in \eqref{pnormestimate}.
\\
Using that $0\leq w_L=(uu_L^{\gamma-1})\leq u^{\gamma}$  on the right-hand side of \eqref{wL} and then letting $L\rightarrow\infty$
on the left-hand side, as a consequence of Fatou's Lemma on the variable $L$, we have
\begin{eqnarray*}
\Biggl(\int_{\Omega}u^{p^*\gamma}dx\Biggl)^{\frac{p}{p^*}}\leq (2\theta_p)^p\gamma^p\Biggl(\int_{\Omega}u^{\gamma\alpha^*}dx\Biggl)^{\frac{p}{\alpha^*}},
\end{eqnarray*}
from which we get that
\begin{eqnarray}\label{gamma}
|u|_{L^{p^*\gamma}(\Omega)}\leq(2\theta_p)^{\frac{1}{\gamma}}\gamma^{\frac{1}{\gamma}}|u|_{L^{\gamma\alpha^*}(\Omega)}.
\end{eqnarray}
\\
Let us define $\sigma:=p^*/\alpha^*.$ When $\gamma=\sigma$ in \eqref{gamma}, since $\gamma\alpha^*=p^*$ we have $u\in L^{p^*\sigma}(\Omega)$ and
\begin{eqnarray}\label{gamma1}
|u|_{L^{p^*\sigma}(\Omega)}\leq (2\theta_p)^{\frac{1}{\sigma}}\sigma^{\frac{1}{\sigma}}
|u|_{L^{p^*}(\Omega)}.
\end{eqnarray}
Now, choosing $\gamma=\sigma^2$ in \eqref{gamma}, since $\gamma\alpha^*=p^*\sigma $ and $p^*\gamma=p^*\sigma^2,$ we obtain
\begin{eqnarray}\label{gamma2}
|u|_{L^{p^*\sigma^2}(\Omega)}\leq (2\theta_p)^{\frac{1}{\sigma^2}}\sigma^{\frac{2}{\sigma^2}}
|u|_{L^{p^*\sigma}(\Omega)},
\end{eqnarray}
by using \eqref{gamma1} and \eqref{gamma2}, we have
\begin{align*}
|u|_{L^{p^*\sigma^2}(\Omega)}\leq(2\theta_p)^{\frac{1}{\sigma^2}+\frac{1}{\sigma}}
\sigma^{\frac{2}{\sigma^2}+{\frac{1}{\sigma}}}
|u|_{L^{p^*}(\Omega)}.
\end{align*}
For $n\geq 1,$ we define $\sigma_n$ inductively so that $\sigma_n=\sigma^n.$ Then, from \eqref{gamma}, it follows that
\begin{eqnarray}\label{gamma-n}
|u|_{L^{p^*\sigma^n}(\Omega)}\leq (2\theta_p)^{\frac{1}{\sigma^n}+...+\frac{1}{\sigma^2}+\frac{1}{\sigma}}
\sigma^{\frac{n}{\sigma^n}+...+\frac{2}{\sigma^2}+{\frac{1}{\sigma}}}
|u|_{L^{p^*}(\Omega)}.
\end{eqnarray}
Note that
$$
\displaystyle\sum_{i=1}^{\infty}\frac{1}{\sigma^i}=\frac{1}{\sigma-1} \ \mbox{ and } \
\displaystyle\sum_{i=1}^{\infty}\frac{i}{\sigma^i}=\frac{\sigma}{(\sigma-1)^2}.
$$
Thus, since $\sigma>1,$ passing to the limit as $n\rightarrow\infty$ in \eqref{gamma-n} we conclude that $u\in L^{\infty}(\Omega)$ and
\begin{eqnarray}\label{est00}
|u|_{L^\infty(\Omega)}\leq (2\theta_p)^{\frac{1}{\sigma-1}}\sigma^{\frac{\sigma}{(\sigma-1)^2}}|u|_{L^{p^*}(\Omega)}.
\end{eqnarray}
Finally, since $\sigma=\frac{N}{N-p}-\frac{q}{p}+1$ and
$1<p<q<1^*<p^*,$ using once more \eqref{p*}, the expression of $\theta_p$ (see \eqref{alphatheta}) and \eqref{est00}
we conclude the proof of the lemma.
\fim

As a consequence of Lemma \ref{bounded} and the
compactness of the embedding $BV(\Omega)\hookrightarrow L^r(\Omega)$, for $ r\in[1,1^*)$ (where $1^* := N/(N-1)$),
it follows that there exists $u_\beta \in BV(\Omega)$ such that, as  $p \to 1^+,$
\begin{equation}
u_{p,\beta} \to u_\beta \quad \mbox{ in } L^r(\Omega)
\label{upconvergence}
\end{equation}
and
\begin{equation}
u_{p,\beta}(x) \to u_\beta(x) \mbox{ a.e. in } \Omega.
\label{upconvergenceae}
\end{equation}

Hence, according to Lemma \ref{regularity1} we have $u_\beta\in L^\infty(\Omega)$ and $u_\beta(x)\geq0$ for almost every $x\in\Omega.$

In what follows, we will prove that $u_\beta$ is a solution of \eqref{P1intro}, in the sense of definition \eqref{sol.u}. Furthermore, we will prove that  $u_\beta\not\equiv0.$

We start with the following result:

\begin{lemma}
\label{lemmaconvrhopbeta}
Let $u_{p,\beta} \in W^{1,p}_0(\Omega),$ $\rho_{p,\beta} \in L^{\frac{q}{q-1}}(\Omega)$ and
$u_\beta\in BV(\Omega)$
satisfying \eqref{Ppweak}, \eqref{rhop} and \eqref{upconvergence}. Then, there exists $\rho_{\beta} \in L^{\frac{q}{q-1}}(\Omega)$, such that
\begin{equation}
\rho_{p,\beta} \rightharpoonup \rho_{\beta} \mbox{ in }  L^{\frac{q}{q-1}}(\Omega),  \mbox{ as p } \to 1^+,
\label{rhopweakconvergence}
\end{equation}
\begin{equation}
\rho_{p,\beta}(x) \to \rho_{\beta}(x)  \mbox{ a.e. in } \Omega,  \mbox{ as p } \to 1^+
\label{rhopae}
\end{equation}
and
\begin{equation}
0 \leq \rho_{\beta}(x) \leq |u_\beta(x)|^{q-1} \mbox{ a.e. in } \Omega.
\label{rhou0}
\end{equation}

Moreover, $\rho_{\beta}$ satisfies, for almost every $x\in \Omega,$
\begin{equation}
\rho_{\beta}(x) \in
\begin{cases}
\{0\}, & \ \mbox{if} \ u_\beta(x)<\beta,\\
 [0, \beta^{q-1}], & \ \mbox{if} \ u_\beta(x)=\beta,\\
 \{u_\beta(x)^{q-1}\}, & \ \mbox{if} \ u_\beta(x)>\beta.
 \end{cases}\label{rhocondition}
\end{equation}

\end{lemma}
\dem
By Lemma \ref{regularity1}, it follows that $(\rho_{p,\beta})_{1 < p < \overline{p}}$ is bounded in $L^{\frac{q}{q-1}}(\Omega)$. Hence, there exists $\rho_{\beta} \in L^{\frac{q}{q-1}}(\Omega)$, such that \eqref{rhopweakconvergence} holds. Moreover, if $E \subset \Omega$ is a measurable set, then
$$
\int_E \rho_{p,\beta} dx = \int_\Omega \rho_{p,\beta} . \chi_E dx \to \int_\Omega \rho_{\beta} . \chi_E dx= \int_E \rho_{\beta} dx, \mbox{ as } p \to 1^+.
$$
Hence, \eqref{rhopae} holds.

Moreover, again by \eqref{rhop} and \eqref{upconvergenceae}, \eqref{rhopae} implies that $\rho_{\beta}$ satisfies \eqref{rhou0}.

Now, we will prove the relation \eqref{rhocondition}
between $\rho_{\beta}$ and $u_\beta.$ We first show that $\rho_{\beta}(x)=0$ for almost every $x\in \{u_\beta<\beta\}$. Note that from \eqref{upconvergenceae} and \eqref{rhou0},
it follows that, as  $p\rightarrow 1^+,$
\begin{equation}\label{rh0}
u_{p,\beta}(x)\rightarrow u_\beta(x) \mbox{ and } \rho_{p,\beta}(x)\rightarrow \rho_{\beta} (x), \mbox{a.e. in $\{u_\beta<\beta\}$}.
\end{equation}
Using \eqref{rh0} we conclude that, for almost every $x\in \{u_\beta<\beta\}$ there exists $p_x>0$ sufficiently small such that $u_{p,\beta}(x)<\beta$ for all $1<p<p_x.$ By the relation \eqref{rhop}, we have $\rho_{p,\beta}(x)=0$ for all $1<p<p_x,$ and thus
by \eqref{rh0} we have $\rho_{\beta}(x)=0.$ Hence, $\rho_{\beta}(x)=0$ for almost every $x\in \{u_\beta<\beta\}.$

By the same arguments as above we can show that $\rho_{\beta}(x)=u_\beta(x)^{q-1}$ for almost every $x\in \{u_\beta>\beta\}.$
Moreover, from \eqref{rhou0}  it follows that $0\leq \rho_{\beta}(x)\leq \beta^{q-1}$ for almost every $x\in \{u_\beta=\beta\}.$
Therefore, we conclude that \eqref{rhocondition} holds true.

\fim

\begin{lemma}\label{z} For each $\beta>0,$ there exists a vector field $\textbf{z}_{\beta} \in L^\infty(\Omega,\mathbb{R}^N)$ such that $\|\textbf{z}_{\beta}\|_\infty \leq 1$ and
\begin{equation}
-\mbox{div}\,  \textbf{z}_{\beta}= \rho_{\beta}, \quad \mbox{in $\mathcal{D}'(\Omega)$,}
\label{equationzu}
\end{equation}
with $\rho_{\beta}$ satisfying \eqref{rhocondition}.
\end{lemma}
\dem

The inequality \eqref{pnormestimate} implies that (see \cite[Proposition 3]{AndreuBallesterCasellesMazon} or \cite[Theorem 3.3]{Mercaldo-Rossi-Segura}) there exists $\textbf{z}_{\beta} \in L^\infty(\Omega,\mathbb{R}^N)$, such that $\|\textbf{z}_{\beta}\|_\infty \leq 1$ and
\begin{equation}
|\nabla u_{p,\beta}|^{p-2}\nabla u_{p,\beta} \rightharpoonup \textbf{z}_{\beta} \mbox{ weakly in } L^r(\Omega,\mathbb{R}^N), \mbox{ as } p \to 1^+,
\label{vectorconvergence}
\end{equation}
for all $1 \leq r < \infty$. In particular, as $p \to 1^+$,
\begin{equation}
|\nabla u_{p,\beta}|^{p-2}\nabla u_{p,\beta} \to \mbox{div}\,  \textbf{z}_{\beta} \quad \mbox{in $\mathcal{D}'(\Omega)$}.
\label{divvectorconvergence}
\end{equation}

Therefore, by using \eqref{Ppweak}, \eqref{rhopweakconvergence} and \eqref{divvectorconvergence} and the Lebesgue dominated convergence theorem, we conclude that
\begin{equation*}
-\mbox{div}\,  \textbf{z}_{\beta} = \rho_{\beta}, \quad \mbox{in $\mathcal{D}'(\Omega)$,}
\end{equation*}
which proves the lemma.
\fim

\begin{lemma} The function $u_\beta$ and the vector field $\textbf{z}_{\beta}$ satisfy the following equality in the sense of measures in $\Omega$,
$$
(\textbf{z}_{\beta},D u_\beta) = |D u_\beta|.
$$
\label{lemmazDu}
\end{lemma}
\dem First of all, since $\|\textbf {z}_{\beta}\|_\infty \leq 1$, it follows that, $(\textbf{z}_{\beta},Du_\beta) \leq |Du_\beta|$ in $\mathcal{M}(\Omega)$.
In fact, for any Borel set  $B$, by \eqref{zDuinequality},
\begin{eqnarray*}
\int_B (\textbf{z}_{\beta},Du_\beta) & \leq & \left| \int_B (\textbf{z}_{\beta},Du_\beta) \right|\\
& \leq & \|\textbf{z}_{\beta}\|_\infty \int_B |Du_\beta|\\
& \leq & \int_B |Du_\beta|.
\end{eqnarray*}
Hence, it is enough to show the opposite inequality, i.e., that for all $\varphi \in C^1_0(\Omega)$, $\varphi \geq 0$,
\begin{equation}
\langle (\textbf{z}_{\beta}, Du_\beta), \varphi \rangle \geq \int_\Omega \varphi |Du_\beta|.
\label{inequalityzDu}
\end{equation}
In order to do so, let us consider $u_{p,\beta} \varphi \in W^{1,p}(\Omega)$ as a test function in \eqref{Pp}. Thus we obtain,
\begin{equation}
\int_\Omega \varphi |\nabla u_{p,\beta}|^p dx + \int_\Omega u_{p,\beta}|\nabla u_{p,\beta}|^{p-2}\nabla u_{p,\beta} \cdot \nabla \varphi dx = \int_\Omega \rho_{p,\beta} \varphi dx.
\label{inequalityzDu1}
\end{equation}
Now we shall calculate the lower limit as $p \to 1^+$ in both sides of \eqref{inequalityzDu1}. Before it, note that, Young's inequality and the lower semicontinuity of the map $v \mapsto \displaystyle\int_\Omega \varphi |Dv|$ with respect to the $L^r(\Omega)$ convergence, imply that
\begin{eqnarray*}
\int_\Omega \varphi |Du_\beta| & \leq & \liminf_{p \to 1^+} \int_\Omega \varphi |\nabla u_{p,\beta}| dx\\
& \leq & \liminf_{p \to 1^+}\left( \frac{1}{p} \int_\Omega \varphi |\nabla u_{p,\beta}|^p dx + \frac{p-1}{p}\int_\Omega \varphi dx \right)\\
& = & \liminf_{p \to 1^+} \int_\Omega \varphi |\nabla u_{p,\beta}|^p dx.
\end{eqnarray*}
Moreover, by \eqref{divvectorconvergence}, it follows that
\begin{equation}
\lim_{p \to1^+}\int_\Omega u_{p,\beta}|\nabla u_{p,\beta}|^{p-2}\nabla u_{p,\beta} \nabla \varphi dx = \int_\Omega u_\beta \textbf{z}_{\beta}\cdot \nabla \varphi dx.
\label{inequalityzDu3}
\end{equation}
Finally, Lebesgue's dominated convergence theorem and \eqref{upconvergence} imply that
\begin{equation}
\lim_{p \to 1^+}\int_\Omega \rho_{p,\beta} \varphi dx = \int_\Omega  \rho_{\beta} \varphi dx.
\label{inequalityzDu4}
\end{equation}
Then, from \eqref{equationzu}, \eqref{inequalityzDu1}, \eqref{inequalityzDu3} and \eqref{inequalityzDu4}, it follows that
\begin{eqnarray*}
\langle (\textbf{z}_{\beta}, Du_\beta), \varphi\rangle &=&  -\int_\Omega \varphi u_\beta \mbox{div}\textbf{z}_{\beta} - \int_\Omega u_\beta \textbf{z}_{\beta}\cdot\nabla \varphi dx\\
& =&  \int_\Omega  \rho_{\beta} u_\beta \varphi dx - \int_\Omega u_\beta \textbf{z}_{\beta}\cdot\nabla \varphi dx\\
& = & \lim_{p \to 1^+}\left(\int_\Omega  \rho_{p,\beta} u_{p,\beta} \varphi dx - \int_\Omega u_{p,\beta}|\nabla u_{p,\beta}|^{p-2}\nabla u_{p,\beta} \cdot\nabla \varphi dx\right)\\
& = &  \liminf_{p \to 1^+}\int_\Omega \varphi|\nabla u_{p,\beta}|^p dx\\
& \geq &  \int_\Omega \varphi |Du_\beta|.
\end{eqnarray*}
Then, \eqref{inequalityzDu} holds and this finishes the proof.
\fim

\begin{lemma}
The function $u_\beta$ satisfies $\left[\textbf{z}_{\beta},\nu \right] \in \mbox{sign}(-u_\beta)$ on $\partial \Omega$.
\label{lemmaboundaryRad}
\end{lemma}
\dem
To check  that $\left[\textbf{z}_{\beta},\nu\right]\in sign(-u_\beta) $
it is enough to show that
\begin{equation}
\int_\Omega (| u_\beta | + u_\beta[\textbf{z}_{\beta},\nu])d\mathcal{H}^{N-1} = 0.
\label{eqboundary1}
\end{equation}
Indeed, since
\begin{eqnarray*}
-u_\beta[\textbf{z}_{\beta},\nu]
&\leq&
\| \textbf{z}_{\beta} \|_{L^{\infty}(\Omega)}| u_\beta | \\
&\leq&
| u_\beta |,
\end{eqnarray*}
the integrand in \eqref{eqboundary1} is nonnegative. Then, \eqref{eqboundary1} holds if and only if $[\textbf{z}_{\beta},\nu](-u_\beta) = |u_\beta|$ $\mathcal{H}^{N-1}$a.e. on $\partial \Omega$.

In order to verify \eqref{eqboundary1}, let us consider $(u_{p,\beta} - \varphi) \in W^{1,p}_{0}(\Omega)$ as test function in \eqref{Pp} with $\varphi \in C^1_{0}(\Omega)$. Then we get
\begin{equation}
\int_\Omega |\nabla u_{p,\beta}|^p dx = \int_\Omega |\nabla u_{p,\beta}|^{p-2}\nabla u_{p,\beta} \nabla \varphi dx + \int_\Omega \rho_{p,\beta} (u_{p,\beta} - \varphi)dx.
\label{eqboundary3}
\end{equation}
From Young's inequality, Green's Formula, \eqref{equationzu}, \eqref{divvectorconvergence}, Lemma \ref{lemmazDu} and \eqref{eqboundary3}, we have that, as $p \to 1^+$,
\begin{eqnarray}
\nonumber p \int_\Omega |\nabla u_{p,\beta}| dx & \leq & \int_\Omega |\nabla u_{p,\beta}|^p dx + (p-1)|\Omega|\\
\nonumber & = & \int_\Omega |\nabla u_{p,\beta}|^{p-2}\nabla u_{p,\beta} \nabla \varphi dx + \int_\Omega \rho_{p,\beta}(u_{p,\beta} - \varphi)dx + (p-1)|\Omega|\\
\nonumber & = & \int_\Omega \textbf{z}_{\beta}\cdot \nabla \varphi dx +  \int_\Omega \rho_{\beta}(u_\beta - \varphi)dx + o_p(1)\\
\nonumber & = & -\int_\Omega \varphi \mbox{div}\textbf{z}_{\beta} -\int_\Omega \rho_{\beta} \varphi dx + \int_\Omega \rho_{\beta} u_\beta dx + o_p(1)\\
\nonumber & = & \int_\Omega \rho_{\beta} u_\beta dx + o_p(1)\\
\nonumber & = & -\int_\Omega u_\beta \mbox{div}\textbf{z}_{\beta} + o_p(1)\\
\nonumber & = & \int_\Omega (\textbf{z}_{\beta}, Du_\beta) - \int_{\partial \Omega} \left[\textbf{z}_{\beta},\nu \right] u_\beta d\mathcal{H}^{N-1} + o_p(1)\\
\label{eqboundary4} & = & \int_\Omega |Du_\beta| - \int_{\partial \Omega} \left[\textbf{z}_{\beta},\nu \right] u_\beta d\mathcal{H}^{N-1} + o_p(1).
\end{eqnarray}
Hence, from \eqref{eqboundary4} and the lower semicontinuity of the norm in $BV(\Omega)$, it follows that
\begin{equation}
\int_{\partial B}\left(|u_\beta| + \left[\textbf{z}_{\beta},\nu \right] u_\beta \right) d\mathcal{H}^{N-1} \leq 0.
\label{eqboundary5}
\end{equation}
But the last inequality imply in \eqref{eqboundary1} and we are done.

\fim

Now, let us prove that the function
$u_\beta \in BV(\Omega)\cap L^\infty(\Omega)$ is a nonnegative and nontrivial solution of \eqref{P1intro}, in the sense of the definition \eqref{sol.u}.

First of all, note that by Lemmas \ref{z}, \ref{lemmazDu} and \ref{lemmaboundaryRad},
$u_\beta\in BV(\Omega),$ $\rho_{\beta}\in L^{\frac{q}{q-1}}(\Omega)$ and $\textbf{z}_{\beta}\in L^{\infty}(\Omega,\mathbb{R}^N)$ satisfy \eqref{sol.u} and \eqref{sol.rho}. Moreover, since $u_{p,\beta}(x)\geq0$
for almost every $x\in \Omega,$
according to \eqref{upconvergenceae} and Lemma \ref{regularity1}, it follows that $u_\beta\in L^{\infty}(\Omega)$
and $u_\beta(x)\geq0$ for almost every $x\in \Omega.$
\\
Now let us show that $u_\beta\not\equiv0.$
Invoking Lemma \ref{geometria} and \eqref{cp}, we have
\begin{equation}\label{ineq-cp}
\alpha + o_p(1) \leq c_{p,\beta}=I_{H,p}(u_{p,\beta})\leq \frac{1}{p}\int_{\Omega}|\nabla u_{p,\beta}|^p dx + o_p(1).
\end{equation}
Hence, since  $\alpha$ is independent of $p$ (see Lemma \ref{geometria}),  \eqref{Ppweak}, \eqref{rhopae}, \eqref{rhou0} and Lebesgue's dominated convergence theorem, imply that
\begin{equation}\label{ineq-alpha}
	\alpha\leq\lim_{p\rightarrow1^+}\frac{1}{p}\int_{\Omega}|\nabla u_{p,\beta}|^p dx=\lim_{p\rightarrow1^+}\frac{1}{p}\int_{\Omega}\rho_{p,\beta}u_{p,\beta} dx=\int_{\Omega}\rho_{\beta}u_\beta dx.
\end{equation}
Thus, combining \eqref{equationzu}, Green's formula (see \eqref{GreenFormula}), Lemma \ref{lemmazDu}, Lemma \ref{lemmaboundaryRad} and \eqref{ineq-alpha}, we deduce that
\begin{eqnarray*}
0& <& \alpha \\
& \leq  & \int_{\Omega}(\textbf{z}_{\beta},Du_\beta) - \int_{\partial\Omega}[\textbf{z}_\beta,\nu]u_\beta\;d\mathcal{H}^{N-1}\\
    &= & \int_{\Omega}|Du_\beta| - \int_{\partial\Omega}[ \textbf{z}_\beta,\nu]u_\beta\;d\mathcal{H}^{N-1}\\
    & = & \int_{\Omega}|Du_\beta| + \int_{\partial\Omega}|u_\beta|\;d\mathcal{H}^{N-1}\\
    & = & \|u_\beta\|,
\end{eqnarray*}
thus $u_\beta\not\equiv0$. Then Theorem \ref{existenceP1} is proved.

\section{Proof of Theorem \ref{theorem2}}
\label{Prooftheorem2}

Now, let us perform a deep analysis of the behavior of $u_{p,\beta}$, as $\beta \to 0^+$.

For each $\beta > 0$, let us define the functional $I_\beta: BV(\Omega) \to \mathbb{R}$, given by
$$
I_\beta (u) = \int_{\Omega}|Du| + \int_{\partial \Omega}|u| d\mathcal{H}^{N-1} - \int_\Omega F_\beta(u) dx.
$$

Note that, since ${\bf z}_{\beta}$ and $u_\beta$ satisfy
\begin{equation}
\label{u0betasolution}
\left\{
\begin{array}{rcl}
- \mbox{div}\, {\textbf z}_{\beta} & = & \rho_{\beta} \quad  \mbox{ in $\mathcal{D}'(\Omega)$,}\\
({\textbf z}_{\beta},Du_\beta) & = & |Du_\beta| \quad \mbox{in $\mathcal{M}(\Omega)$},\\
{[{\textbf{z}_{\beta}},\nu]} & \in & \mbox{$sign(-u_\beta)$ \quad $\mathcal{H}^{N-1}$-a.e. on $\partial \Omega,$}
\end{array}
 \right.
\end{equation}
by taking $u_\beta$ as test function in \eqref{u0betasolution} and using Green Formula, \eqref{Pp} and \eqref{rhopweakconvergence}, it follows that
\begin{eqnarray}
\nonumber \|u_\beta\| & = & \int_\Omega |Du_\beta| + \int_{\partial \Omega}|u_\beta| d\mathcal{H}^{N-1}\\
\nonumber & = & -\int_\Omega u_\beta\mbox{div}\, 
{\textbf{z}_{\beta}}\\
\label{convu01} & = & \int_\Omega u_\beta\rho_{\beta} dx\\
\nonumber & = & \int_\Omega u_{p,\beta}\rho_{p,\beta} dx + o_p(1)\\
\nonumber & = & \int_\Omega |\nabla u_{p,\beta}|^p dx + o_p(1).
\end{eqnarray}

Moreover, from \eqref{upconvergence} and \eqref{upconvergenceae}, it follows that
\begin{equation}
\label{convu02}
\int_\Omega F_\beta(u_\beta) dx = \int_\Omega F_\beta(u_{p,\beta}) dx + o_p(1).
\end{equation}
Hence, from \eqref{convu01} and \eqref{convu02}, we have that
\begin{equation}
\label{convu03}
I_{\beta}(u_\beta) = I_{p,\beta}(u_{p,\beta}) + o_p(1).
\end{equation}

Since we are interested in the behavior of $u_\beta$, as $\beta \to 0^+$, let us assume from now on that $0 < \beta < \beta_0$.

\begin{lemma}
The family $(u_\beta)_{0 < \beta < \beta_0}$ is bounded in $BV(\Omega)$.
\end{lemma}
\dem
First of all, let us prove that, if $0 < \beta_1 < \beta_2 < \beta_0$, then
\begin{equation}
\label{ubetalim0}
I_{\beta_1}(u_{\beta_1}) \leq I_{\beta_2}(u_{\beta_2}).
\end{equation}
In order to do so, let us prove that, for $p > 1$ fixed,
\begin{equation}
\label{ubetalim1}
I_{p,\beta_1}(u_{p,\beta_1}) < I_{p,\beta_2}(u_{p,\beta_2}).
\end{equation}
Note that, for $u \in W^{1,p}_0(\Omega)$, since $F_{\beta_1}(u) \geq F_{\beta_2}(u)$ a.e. in $\Omega$, it follows that
\begin{equation}
\label{ubetalim2}
I_{p,\beta_1}(u) \leq I_{p,\beta_2}(u).
\end{equation}
Moreover, let us assume that the function $e$ in Lemma \ref{geometria}, is $e(\beta_0)$, i.e., that satisfies $I_{p,\beta_0}(e) < p/(p-1)$. Hence, from \eqref{ubetalim2}, we have that
$$
I_{p,\beta}(e) \leq I_{p,\beta_0}(e) < \frac{p}{p-1},
$$
for all $0 < \beta < \beta_0$.  Hence, in the definition of $c_{p,\beta}$, for $0 < \beta < \beta_0$, we can assume without loss of generality that $e = e(\beta_0)$ and then the class of paths $\Gamma$ does not depend on $\beta$. Then, from \eqref{ubetalim2}, it follows that
\begin{eqnarray*}
I_{p,\beta_1}(u_{p,\beta_1}) & = & c_{p,\beta_1}\\
 & = & \inf_{\gamma \in \Gamma}\sup_{t \in [0,1]} I_{p,\beta_1}(\gamma(t))\\
& \leq & \inf_{\gamma \in \Gamma}\sup_{t \in [0,1]} I_{p,\beta_2}(\gamma(t))\\
& = & c_{p,\beta_2}\\
& = & I_{p,\beta_2}(u_{p,\beta_2}).
\end{eqnarray*}
This, in turn, proves \eqref{ubetalim1}.

Hence, from \eqref{convu03}, passing the limit as $p \to 1^+$ in \eqref{ubetalim1}, we have that \eqref{ubetalim0} holds.

Then, for all $0 < \beta < \beta_0$,
$$
I_\beta(u_\beta) \leq I_{\beta_0}(u_{\beta_0}) =: C.
$$

Note that, by using $u_\beta$ as test function in \eqref{u0betasolution}, from Green Formula, we have that
\begin{eqnarray}
\nonumber \|u_\beta\| & = & \int_\Omega |Du_\beta| + \int_{\partial \Omega}|u_\beta| d\mathcal{H}^{N-1}\\
\nonumber & = & \int_\Omega ({\bf z}_\beta,Du_\beta) + \int_{\partial \Omega}[{\bf z}_\beta,\nu]u_\beta d\mathcal{H}^{N-1}\\
\label{ubetaconv1}& = & -\int_\Omega u_\beta \mbox{div}\, {\bf z}_\beta dx\\
\nonumber & = & \int_\Omega u_\beta \rho_\beta dx.
\end{eqnarray}

Then, from \eqref{rhocondition}, \eqref{ubetaconv1} and the definition of $F_\beta$, we have that

\begin{eqnarray*}
I_\beta(u_\beta)  & = & I_\beta(u_\beta) - \frac{1}{q}\left(\|u_\beta\| - \int_\Omega u_\beta \rho_\beta dx\right)\\
& = & \left(1 - \frac{1}{q}\right)\|u_\beta\| + \int_\Omega \left(\frac{1}{q} u_\beta \rho_\beta - F_\beta(u_\beta)\right) dx\\
& \geq &  \left(1 - \frac{1}{q}\right)\|u_\beta\| + \frac{\beta^q}{q}|\{u_\beta > \beta\}|\\
& \geq & \left(1 - \frac{1}{q}\right)\|u_\beta\|.
\end{eqnarray*}
Hence, since $(I_\beta(u_{\beta}))_{0 < \beta < \beta_0}$ is bounded, it follows from the last inequality that $(u_\beta)_{0 < \beta < \beta_0}$ is also bounded.

\fim

From the last result, there exists $u_0 \in BV(\Omega)$ such that, for all $r \in [1,1^*)$,
\begin{equation}
u_{\beta} \to u_0 \quad \mbox{ in } L^r(\Omega)
\label{ubetaconvergence}
\end{equation}
and
\begin{equation}
u_{\beta}(x) \to u_0(x) \mbox{ a.e. in } \Omega.
\label{ubetaconvergenceae}
\end{equation}
Moreover, note that the boundedness on $(u_\beta)_{0 < \beta < \beta_0}$ and \eqref{rhocondition} implies also that $(\rho_\beta)_{0 < \beta < \beta_0}$ is bounded in $L^\frac{q}{q-1}(\Omega)$. Then, as in Lemma \ref{lemmaconvrhopbeta}, it is possible to show that there exists $\rho_0 \in L^\frac{q}{q-1}(\Omega)$, such that
\begin{equation}
\rho_{\beta} \rightharpoonup \rho_0 \mbox{ in }  L^{\frac{q}{q-1}}(\Omega),  \mbox{ as $\beta \to 0^+ $},
\label{rhobetaweakconvergence}
\end{equation}
\begin{equation}
\rho_{\beta}(x) \to \rho_0(x)  \mbox{ a.e. in } \Omega,  \mbox{ as $\beta \to 0^+$}
\label{rhobetaae}
\end{equation}
and
\begin{equation}
0 \leq \rho_0(x) \leq |u_0(x)|^{q-1} \mbox{ a.e. in } \Omega.
\label{rho0}
\end{equation}

Now, let us deal with the family of vector fields $\left({\bf z}_\beta\right)_{0 < \beta < \beta_0}$. Note that, since $\|{\bf z}_\beta\|_\infty \leq 1$ for all $\beta \in (0,\beta_0)$, then there exists ${\bf z}_0 \in L^\infty(\Omega,\mathbb{R}^N)$, such that
\begin{equation}
{\bf z}_\beta \overset{\ast}{\rightharpoonup} {\bf z}_0 \quad \mbox{in $L^\infty(\Omega,\mathbb{R}^N)$}.
\end{equation}
This, on the other hand, implies that ${\bf z}_\beta \rightharpoonup {\bf z}_0$ in $L^1(\Omega,\mathbb{R}^N)$, i.e., for all $\psi \in L^\infty(\Omega,\mathbb{R}^N)$,
\begin{equation}
\int_\Omega {\bf z}_\beta \cdot \psi dx \to \int_\Omega {\bf z}_0 \cdot \psi dx, \quad \mbox{as $\beta \to 0^+$.}
\label{z_betaconvergence}
\end{equation}
For every $\varphi \in C^\infty_c(\Omega)$, since $\nabla \varphi \in L^\infty(\Omega,\mathbb{R}^N)$, by \eqref{z_betaconvergence}, we have that
$$
\int_\Omega {\bf z}_\beta\cdot \nabla\phi dx \to \int_\Omega {\bf z}_0 \cdot \nabla \phi dx, \quad \mbox{as $\beta \to 0^+$},
$$
from where it follows that
\begin{equation}
\mbox{div}\, {\bf z}_\beta \to \mbox{div}\, {\bf z}_0, \quad \mbox{in $\mathcal{D}'(\Omega)$.}
\label{z_betaconvergence2}
\end{equation}

Hence, from \eqref{rhobetaweakconvergence} and \eqref{z_betaconvergence2}, we have that
\begin{equation}
\label{u0solution1}
-\mbox{div} \, {\bf z}_0 = \rho_0 \quad \mbox{in $\mathcal{D}'(\Omega)$}.
\end{equation}

\begin{lemma}
The function $u_0$ and the vector field ${\bf z}_0$ satisfy the following equality,
$$
({\bf z}_0,D u_0) = |D u_0| \quad \mbox{in $\mathcal{M}(\Omega)$}.
$$
\label{lemmaz0Du0}
\end{lemma}
\dem
First of all, note that, from \eqref{zDuinequality},
$$
({\bf z}_0,Du_0) \leq |Du_0| \quad \mbox{in $\mathcal{M}(\Omega)$.}
$$
For the inverse inequality, let $\varphi \in \mathcal{D}(\Omega)$, $\varphi \geq 0$. In \eqref{u0betasolution}, let us take $\varphi u_\beta$ as test function in \eqref{u0betasolution}. Then,
$$
\int_\Omega \varphi ({\bf z}_\beta,Du_\beta) = \int_\Omega \varphi u_\beta \rho_\beta dx - \int_\Omega u_\beta {\bf z}_\beta\cdot \nabla \varphi dx.
$$
Taking into account that $({\bf z}_\beta,Du_\beta) = |Du_\beta|$ in $\mathcal{M}(\Omega)$,
$$
\int_\Omega \varphi |Du_\beta| = \int_\Omega \varphi u_\beta \rho_\beta dx - \int_\Omega u_\beta {\bf z}_\beta\cdot \nabla \varphi dx.
$$
Taking the $\liminf$ as $\beta \to 0^+$, from the lower semicontinuity of the norm in $BV(\Omega)$ with respect to the $L^r$ convergence, \eqref{ubetaconvergence}, \eqref{rhobetaweakconvergence}  and \eqref{z_betaconvergence2}, it follows that
\begin{eqnarray*}
\int_\Omega \varphi |Du_0| & \leq & \liminf_{\beta \to 0^+} \left(-\int_\Omega \varphi u_\beta \mbox{div}\, {\bf z}_\beta dx - \int_\Omega u_\beta {\bf z}_\beta\cdot \nabla \varphi dx\right)\\
& = & -\int_\Omega \varphi u_0 \mbox{div}\, {\bf z}_0 dx - \int_\Omega u_0 {\bf z}_0\cdot \nabla \varphi dx\\
& = & \int_\Omega \varphi ({\bf z}_0, Du_0).
\end{eqnarray*}
This, in turn, proves that $|Du_0| \leq ({\bf z}_0,Du_0)$ and this finishes the proof.
\fim

\begin{lemma}
The function $u_0$ satisfies $\left[\textbf{z}_0,\nu \right] \in \mbox{sign}(-u_0)$ on $\partial \Omega$.
\label{lemmaboundaryu0}
\end{lemma}
\dem
As in Lemma \ref{lemmaboundaryRad}, it is enough to show that
\begin{equation}
\int_\Omega (| u_0 | + u_0[\textbf{z}_0,\nu])d\mathcal{H}^{N-1} = 0.
\label{eqboundary1}
\end{equation}

In order to verify \eqref{eqboundary1}, let us consider $(u_\beta - \varphi) \in BV(\Omega)\cap L^\infty(\Omega)$ as test function in \eqref{u0betasolution}, where $\varphi \in \mathcal{D}(\Omega)$. Then, from \eqref{GreenFormula} and \eqref{u0betasolution}, we get
\begin{eqnarray}
\nonumber \int_\Omega |Du_\beta| + \int_{\partial \Omega} |u_\beta| d\mathcal{H}^{N-1} & = & \int_\Omega ({\bf z}_\beta,Du_\beta) - \int_{\partial \Omega} u_\beta [{\bf z}_\beta,\nu] d\mathcal{H}^{N-1}\\
\label{eqboundary2} & = & -\int_\Omega u_\beta \mbox{div}\, {\bf z}_\beta\\
\nonumber & = & \int_\Omega \varphi \mbox{div}\, {\bf z}_\beta + \int_\Omega u_\beta \rho_\beta dx - \int_\Omega \varphi \rho_\beta\\
\nonumber & = & \int_\Omega u_\beta \rho_\beta dx.
\end{eqnarray}
Then, calculating the $\liminf$ in \eqref{eqboundary2}, from the lower semicontinuity of the norm in $BV(\Omega)$, \eqref{u0solution1} and Lemma \ref{lemmaz0Du0}, we have that

\begin{eqnarray*}
\int_\Omega |Du_0| \int_{\partial \Omega}|u_0| d\mathcal{H}^{N-1} & \leq & \int_\Omega u_0 \rho_0 dx\\
& = & -\int_\Omega u_0 \mbox{div}\, {\bf z}_0 \\
& = & \int_\Omega ({\bf z}_0,Du_0) - \int_{\partial \Omega} u_0 [{\bf z}_0,\nu] \mathcal{H}^{N-1} \\
& = & \int_\Omega |Du_0| - \int_{\partial \Omega} u_0 [{\bf z}_0,\nu] \mathcal{H}^{N-1}.
\end{eqnarray*}
From the last inequality, it follows that
$$
| u_0 | + u_0[\textbf{z}_0,\nu] \leq 0 \quad \mbox{$\mathcal{H}^{N-1}-$a.e. on $\partial \Omega$.}
$$
Since the inverse inequality is trivial, it follows that \eqref{eqboundary1} holds.
\fim

Then, from \eqref{u0solution1} and Lemmas \ref{lemmaz0Du0} and \ref{lemmaboundaryu0}, it follows that $u_0$ is a solution of \eqref{P2intro}.

Now, in order to end up the proof of Theorem \ref{theorem2}, let us show that there exist constants $\mu, \beta_{0} > 0$, such that
\begin{equation}
\label{muubeta}
|\{x\in \Omega: u_\beta(x)>\beta\}|\geq
\mu, \ \mbox{ for all } \ \beta\in(0, \beta_{0}),
\end{equation}

From \eqref{ineq-cp} it
follows that
$$
0 < \alpha + o_p(1) \leq c_{p,\beta} \leq \frac{1}{p}\int_{\Omega}\rho_{p,\beta} u_{p,\beta} dx + o_p(1),
$$
where $\alpha$ is independent of $\beta$ and $p\in (1, \overline{p})$. Since $\rho_{p,\beta}$ verifies \eqref{rhop},
\begin{equation}\label{ineq}
\alpha\leq \frac{\beta^{q}}{p}|\Omega|+\frac{1}{p}\int_{\{u_{p,\beta}>\beta\}}u_{p,\beta}^{q} dx + o_p(1),
\end{equation}
for all $\beta>0$ and $p\in (1, \overline{p})$. To conclude the proof, it is enough to prove that
\begin{equation}\label{conver}
 \limsup\limits_{p\to 1^{+}}\int_{\{u_{p,\beta}>\beta\}}u_{p,\beta}^{q} dx\leq \int_{\{u_\beta\geq \beta\}}u_\beta^{q} dx.
\end{equation}
In fact, if \eqref{conver} holds true, then from \eqref{upconvergence}, passing to the upper limit as $p\to 1^{+}$ in \eqref{ineq}, we get
\begin{equation}\label{ineq0}
\alpha\leq \beta^{q}|\Omega|+\int_{\{u_\beta\geq \beta\}}u_\beta^{q} dx\leq 2\beta^{q}|\Omega|+\int_{\{u_\beta> \beta\}}u_\beta^{q} dx,
\end{equation}
for all $\beta>0$. Now, suppose by contradiction that there exists a subsequence $\beta_{n}\to 0$ such that
\begin{equation}\label{conv}
|\{u_\beta>\beta_{n}\}|\to 0, \ \mbox{as \ $\beta_{n}\to 0$}.
\end{equation}
Since, by \eqref{ineq0}, we have
$$
\alpha\leq 2\beta_{n}^{q}|\Omega|+\int_{\Omega}u_\beta^{q}\chi_{\{u_\beta>\beta_{n}\}} dx,
$$
it follows from H\"older's inequality that
\begin{equation}\label{isso}
\alpha\leq 2\beta_{n}^{q}|\Omega|+\left(\int_{\Omega}u_\beta^{r}dx\right)^{\frac{q}{r}}|\{u_\beta>\beta_{n}\}|^{\frac{r-q}{r}},
\end{equation}
for some $q<r<N/(N-1)$. Then, \eqref{conv} and \eqref{isso} would lead us to a contradiction.

Hence, to conclude the proof, it remains us to show \eqref{conver}. For this purpose, observe that
\begin{equation}\label{acaba1}
\int_{\{u_{p,\beta}>\beta\}}u_{p,\beta}^{q} dx=\int_{\Omega}u_{p,\beta}^{q} \chi_{\{u_{p,\beta}>\beta\}}dx\leq \int_{\Omega}u_{p,\beta}^{q} \chi_{\{u_{p,\beta}>\beta\}\cap \{u_\beta<\beta\}}dx+\int_{\Omega}u_{p,\beta}^{q} \chi_{\{u_\beta\geq \beta\}}dx.
\end{equation}
Moreover, since
\begin{equation}\label{acaba2}
\chi_{\{u_{p,\beta}>\beta\}\cap \{u_\beta<\beta\}}(x)\to 0 \ \mbox{a.e. in $\Omega$, \ as \ $p\to 1^{+}$},
\end{equation}
it follows from \eqref{upconvergence}, \eqref{acaba1}, \eqref{acaba2} and Lebesgue's dominated convergence theorem, that \eqref{conver} holds true.

To conclude that $u_0$ is nontrivial, note that by \eqref{ineq0}, $u_0\geq0$ and Lebesgue's dominated convergence theorem,
\begin{equation*}
0<\alpha\leq \int_{\Omega}u_0^{q},
\end{equation*}
then, $u_0\not\equiv0.$ Finally, since $\rho_\beta(x)\in[\underline{f}_\beta(u_\beta(x)),\overline{f}_\beta(u_\beta(x))],$ by \eqref{inclu}, \eqref{ubetaconvergence} and \eqref{rhobetaae},
 we conclude that $\rho_0(x)=u_0(x)^{q-1}$ a.e. in $\Omega$ and therefore $u_0$ is solution of the continuous problem \eqref{P2intro}.

\noindent{\bf Acknowledgments:} Marcos T.O. Pimenta is partially supported by FAPESP 2021/04158-4, CNPq 303788/2018-6 and FAPDF, Brazil.  J. R. Santos J\'unior is partially supported by FAPESP 2021/10791-1 and CNPq 313766/2021-5, Brazil.

\end{document}